\documentclass{article}
\usepackage{ijcai18}

\usepackage{times}
\usepackage{xcolor}
\usepackage{soul}
\usepackage[utf8]{inputenc}
\usepackage[small]{caption}

\usepackage{amssymb,amsmath,amsfonts,graphicx}
\usepackage{hyperref}
\def\X{{\mathbb X}}
\newcommand{\lpr}{{\underline{P}}}

\newtheorem{exam}{Example}
\newtheorem{proposition}{Proposition}
\newtheorem{definition}{Definition}

\title{Statistical inference with belief functions: A survey
}

\author{Fabio Cuzzolin\\ 
Institute for AI, Data Analysis and Systems (AIDAS)\\ School of Engineering, Computing \& Mathematics, Oxford Brookes University  \\
fabio.cuzzolin@brookes.ac.uk}

\begin{document}

\maketitle

\begin{abstract}
Belief functions are a powerful and popular framework for the mathematical characterisation of uncertainty, in particular in situations in which lack of data renders learning a probability distribution for the problem impractical. The first step in a reasoning chain based on belief functions is inference: how to learn a belief measure from the available data. In this survey we focus, in particular, on making inference from statistical data, and review the most significant contributions in the area.
\end{abstract}

\section{Introduction}

Inference is the first step in any estimation or decision problem. In the context of Dempster-Shafer or belief functions theory \cite{Shafer76,cuzzolin2001geometric,cuzzolin2004geometry,cuzzolin2008geometric,cuzzolin14lap,cuzzolin2020geometry,cuzzolin2020geometry-dempster,cuzzolin2014belief}, one of the most popular mathematical frameworks for uncertainty quantification \cite{cuzzolin2021big,cuzzolin2024uncertainty}, inference means {constructing a belief function from the available evidence} \cite{Klopotek1996identification}.\\
Contrarily to Kolmogorov's additive probability measures, belief functions can represent a wide range of `uncertain' data, from classical statistical samples \cite{Seidenfeld78} to qualitative expert judgements, sometimes expressed in terms of mere preferences. 
A number of different approaches in this sense has been proposed: a very general exposition by can be found in \cite{chateauneuf00ambiguity}. Another reference document on inference with belief function is \cite{Ferson03pboxes}, which summarises a variety of the most useful and commonly applied methods for obtaining belief measures, and their mathematical kin probability boxes.

As far as inference from \emph{statistical data}, in particular, is concerned, a number of major streams can be identified: Shafer's approach based on the traditional likelihood function, Wasserman's proposal based on robust Bayesian inference, fiducial inference (including Dempster's proposal based on an auxiliary variable, and weak belief), frequentist methods (in particular Walley and Fine' early work and more recent methods based on confidence intervals and structures). 
\\
The problem can be posed as follows. Consider a \emph{parametric model}, i.e., a family of conditional probability distributions $f(x|\theta)$ of the data $x\in \X$ given a parameter $\theta\in \Theta$:
\begin{equation} \label{eq:parametric-model}
\Big \{ f(x | \theta), x\in \X, \theta\in \Theta \Big \}, 
\end{equation}
where $\X$ is the observation space and $\Theta$ the parameter space. Having observed $x$, how do we quantify the uncertainty about $\theta$, without having to specify a prior distribution?


\section{Statistical inference}

Recall that a \emph{probability measure} over a 
$\sigma$-algebra $\mathcal{F}(\Omega) \subset 2^{\Omega}$, associated with a sample space $\Omega$, is a function $P:\mathcal{F}(\Omega) \rightarrow [0,1]$ such that: $P(\emptyset)=0$; $P(\Omega)=1$; if $A\cap B = \emptyset,\; A,B\in \mathcal{F}(\Omega)$ then $P(A\cup B) = P(A) + P(B)$ (\emph{additivity}).
A sample space $\Omega$ together with a $\sigma$-algebra $\mathcal{F}(\Omega)$ of its subsets and a probability measure $P$ on $\mathcal{F}(\Omega)$ forms a \emph{probability space}, namely the triplet: $(\Omega,\mathcal{F}(\Omega),P)$.
A \emph{random variable} is a function $X$ from a sample space $\Omega$ (endowed with a probability space) to a measurable space $E$ (e.g., $E=\mathbb{R}$).
Statistical inference is then problem of estimating a probability measure given the available sample data.

\textbf{Maximum likelihood estimation} (MLE) \cite{Fisher309} is based on the \emph{likelihood principle}: in a sample, all of the evidence relevant to model parameters is contained in the likelihood function.
Given a \emph{parametric model} (\ref{eq:parametric-model}), the maximum likelihood estimate of $\theta$ is defined as:
${\hat {\theta }}_{\mathrm {MLE} } \subseteq \{ \arg\max_{\theta \in \Theta} \mathcal{L}(\theta \,;\,x_{1},\ldots ,x_{n}) \},$
where the likelihood of the parameter given the observed data $X = \{ x_{1},\ldots ,x_{n} \}$ is:
${\mathcal {L}}(\theta \,;\,x_{1},\ldots ,x_{n}) \doteq f(x_{1},x_{2},\ldots ,x_{n} | \theta ).$

In \textbf{Bayesian inference}, instead, the \emph{prior distribution} is the distribution of the parameter(s) before any data is observed, i.e., $p(\theta | \alpha)$, a function of a vector of hyperparameters $\alpha$.
The \emph{marginal likelihood} (or `evidence') is the distribution of the observed data marginalised over the parameter(s), namely:
$p( {X} | \alpha )=\int _{\theta }p( {X} | \theta )p(\theta | \alpha )\operatorname {d} \!\theta.$
The \emph{posterior distribution} is then the distribution of the parameter(s) after taking into account the observed data, as determined by Bayes' rule:
$
p(\theta |  {X} ,\alpha )={\frac {p( {X} | \theta )p(\theta | \alpha )}{p( {X} | \alpha )}}\propto p( {X} | \theta )p(\theta | \alpha ).
$

\textbf{Frequentist inference} (supported by those who believe probability measures to be limits of relative frequencies) deals with the size of the sample via `confidence intervals'.
\\
Let $X$ be a sample from a probability $P(.|\theta,\phi)$ where $\theta$ is the parameter to be estimated and $\phi$ a nuisance parameter. A \emph{confidence interval} for the parameter $\theta$, with confidence level $\gamma$, is an interval $[u(X),v(X)]$ determined by the pair of random variables $u(X)$ and $v(X)$, with the property:
\begin{equation} \label{eq:confidence-intervals}
{\Pr }(u(X)<\theta <v(X) | \theta ,\phi)=\gamma \quad \forall (\theta ,\phi ). 
\end{equation}
Confidence intervals are a form of interval estimate. Their correct interpretation is about `sampling samples': if we keep extracting new sample sets, 95\% (say) of the time the confidence interval (which will differ for every new sample set) will cover the true value of the parameter. 
One cannot claim, instead, that a specific confidence interval is such that it contains the
value of the parameter with 95\% probability.

\textbf{Fiducial inference} was introduced in \cite{fisher1935fiducial}. 
\\
Let $\theta\in\Theta$ be the parameter of interest, $X$ the observed sample or a sufficient statistic\footnote{A statistic $t =T(X)$ is a function of the observed sample $X$. The statistic is sufficient if $P(X|t,\theta) = P(X|t)$.}, and $u \in U$ an auxiliary variable (termed \emph{pivotal} quantity) such that an \emph{a-equation} holds:
\begin{equation} \label{eq:a-equation}
X = a(\theta,u).
\end{equation}
The crucial assumption underlying the fiducial argument is that each one of $(X,\theta, u)$ is uniquely determined by the a-equation (\ref{eq:a-equation}) given the other two.
The pivotal quantity $u$ is assumed to have an a priori distribution $\mu$, independent of $\theta$.
Prior to the experiment, $X$ has a sampling distribution that depends on $\theta$; after the experiment, however, $X$ is no longer a random variable (as it is measured).
If we `\emph{continue to believe}' that $u$ is distributed according to $\mu$ even after $X$ is observed, we can derive a \emph{fiducial distribution} for $\Theta$.

\begin{exam}\hspace{-1.5mm}\emph{\textbf{: fiducial inference}.} \label{exa:fiducial}
As an example, consider the problem of
estimating the unknown mean of a Gaussian $N(0, 1)$ population based on a single observation $X$.
The a-equation is in this case: $X = \theta + \Psi^{-1}(u)$, where $\Psi(.)$ is the \emph{cumulative distribution function} (CDF)\footnote{The \emph{cumulative distribution function} (CDF) of a random variable $X$ is defined as $F(x) \doteq P(X \leq x)$.} of the $N(0, 1)$ distribution.
Assume the pivotal quantity has prior distribution $\mu = \mathcal{U}(0,1)$.
For a fixed $\theta$, the events $\{\theta \leq \bar{\theta}\}$ and $\{u \geq \Psi(X - \bar{\theta})\}$ are the same. Hence, their probabilities need to be the same.
If we `continue to believe', then, the \emph{fiducial probability} of $\{ \theta \leq \bar{\theta} \}$ is $\Psi(\bar{\theta} - X)$. In other words, the \emph{fiducial distribution of $\theta$ given $X$} is:
$\theta \sim N(X,1)$, a posterior for $\theta$ obtained without requiring a prior for it.
\end{exam}

\section{Belief functions}

Most major approaches to belief function (BF) inference turn out to be generalisations of the above four main strands in statistical inference: inference based on likelihood, Bayesian, fiducial and frequentist inference. 

\textbf{Random set definition}.
Let us denote by $\Omega$ and $\Theta$ the sets of possible answers to two different but related problems $Q_1$ and $Q_2$, respectively. We are given a probability measure $P$ on $\Omega$, and we want to derive a `degree of belief', denoted by $Bel(A)$, that $A\subset \Theta$ contains the correct response to $Q_2$. If we call $\Gamma(\omega)$ the subset of answers to $Q_2$ compatible with $\omega\in\Omega$, each element $\omega$ tells us that the answer to $Q_2$ is somewhere in $A$ whenever $\Gamma(\omega) \subset A$. The \emph{degree of belief} $Bel(A)$ of an event $A\subset\Theta$ is then the total probability \cite{zhou2017total} (in $\Omega$) of all the answers $\omega$ to $Q_1$ that satisfy the above condition, namely \cite{Dempster67}:
\[
Bel(A) = P(\{ \omega | \Gamma(\omega) \subset A \}) = \sum_{\omega\in\Omega
: \Gamma(\omega)\subset A} P(\{\omega\}).
\]
The map $\Gamma : \Omega \rightarrow 2^{\Theta}$ (where $2^\Theta \doteq \{A \subseteq \Theta\}$ denotes the collection of subsets of $\Theta$) is
called a \emph{multivalued mapping} or \emph{random set} from $\Omega$ to $\Theta$. Such a mapping $\Gamma$, together with a probability measure $P$ on $\Omega$, induces a \emph{belief function} on $2^\Theta$.

\textbf{Belief and plausibility measures}.
A \emph{basic probability assignment} (BPA) \cite{Shafer76,cuzzolin08pricai-moebius,cuzzolin10ida} over a finite domain $\Theta$ is a set function 
$m : 2^\Theta\rightarrow[0,1]$ defined on the collection of all subsets of $\Theta$ s.t.:
$m(\emptyset)=0, \; \sum_{A\subset\Theta} m(A)=1.$
The `mass' $m(A)$ assigned to $A$ is nothing but the probability $P(\{\omega\})$ of $\omega \in \Omega : \Gamma(\omega) = A$.
\\
Non-zero mass subsets of $\Theta$ are called \emph{focal elements} of $m$.
\begin{definition} \label{def:bel2}
The \emph{belief function} (BF) associated with a basic probability assignment $m : 2^\Theta\rightarrow[0,1]$ is the set function $Bel : 2^\Theta\rightarrow[0,1]$ defined as:
$Bel(A) = \sum_{B\subseteq A} m(B).$
\end{definition}
As shown in \cite{Shafer76}, belief functions can also be defined axiomatically, without a random set interpretation.
\\
The \emph{plausibility function} $Pl : 2^\Theta \rightarrow [0,1]$ conveys the same information as $Bel$, and is defined as:
$Pl(A) \doteq \sum_{B\cap A\neq \emptyset} m(B) \geq Bel(A).$
The \emph{contour function} associated with $Bel$ is then simply $pl: \Theta \rightarrow [0,1]$, $pl(x) \doteq Pl(\{x\})$.

Classical probability measures on finite sets $\Theta$ are a special case of belief functions (those whose focal elements are singletons only), termed \emph{Bayesian belief functions}.
A belief function is said to be \emph{consonant} if its focal elements $A_1,...,A_m$ are nested: $A_1 \subset A_2 \subset \cdots \subset A_m$.

\textbf{Combination and conditioning}.
The \emph{orthogonal sum} or \emph{Dempster's combination} $Bel_1 \oplus Bel_2 : 2^\Theta \rightarrow [0,1]$ of two belief functions $Bel_1 : 2^\Theta \rightarrow [0,1]$, $Bel_2 : 2^\Theta \rightarrow [0,1]$ defined on the same domain $\Theta$ is the unique BF on $\Theta$ with as focal elements all the {non-empty} intersections of focal elements of $Bel_1$ and $Bel_2$, and basic probability assignment:
$m_{\oplus}(A) = \frac{m_\cap(A)} {1- m_\cap(\emptyset)},$
where $m_i$ denotes the BPA of the input BF $Bel_i$, and:
$m_\cap(A) = \sum_{B \cap C = A} m_1(B) m_2(C).$
Dempster's rule naturally induces a conditioning operator.
Given a conditioning event $A \subset \Theta$, the `categorical' belief function $Bel_A$ such that $m(A)=1$ is combined via Dempster's rule with the a-priori BF $Bel$. The resulting belief function $Bel \oplus Bel_A$ is the conditional BF given $A$ \emph{a la Dempster}, denoted by $Bel_\oplus(A|B)$.
Alternative combination rules have since been proposed \cite{Kramosil02probabilistic-analysis}.

\textbf{Credal interpretation}.
Belief functions admit the following order relation:
$Bel \leq Bel' \equiv Bel(A) \leq Bel'(A)$ $\forall A \subset \Theta$,
called \emph{weak inclusion}. A probability distribution $P$ which weakly includes a belief function $Bel$, $P(A)\geq Bel(A)$ $\forall A \subset \Theta$, is said to be \emph{consistent} with it \cite{kyburg87bayesian}. Each BF thus uniquely identifies a (convex) set of consistent probabilities:
$\mathcal{P}[Bel] = \{ P \in \mathcal{P} : P(A) \geq Bel(A) \}$,
or \emph{credal set} \cite{levi80book}
(where $\mathcal{P}$ is the set of all probabilities one can define on $\Theta$), of which it is its lower envelope. 
Accordingly, BF theory is seen by some authors as a special case of robust statistics (cfr. Section \ref{sec:robust-bayesian-inference}). 
\\
The fact that BFs are mathematically equivalent to lower envelopes $\underline{P}$ of sets of probability measures explains the terminology \emph{lower probability} sometimes applied to these objects: $Bel(A) = \underline{P}(A)$. The corresponding plausibility measure is the \emph{upper probability} of an event $A$: $Pl(A) = \overline{P}(A)$.

The \textbf{least commitment principle} states that, when several belief functions are compatible with a set of constraints, the least informative (according to some informational ordering) should be selected, if such a belief function exists.
For instance, weak inclusion can be adopted as ordering relation.

\textbf{Multivariate analysis}.
Let $\Theta_X$ and $\Theta_Y$ be two sample spaces associated with variables $X,Y$, and let $m^{XY}$ be a mass function on $\Theta_{XY} = \Theta_X \times \Theta_Y$.
The latter can be expressed in the coarser domain $\Theta_X$ by transferring each mass $m^{XY}(A)$ to the {projection} $A\downarrow \Theta_X$ of $A$ on $\Theta_X$.
We obtain a \emph{marginal} mass function on $\Theta_X$, denoted by: $
m^{XY}_{\downarrow X}(B) \doteq \sum_{\{A \subseteq \Theta_{XY}, A\downarrow \Theta_X=B\}}  m^{XY}(A), \forall B \subseteq \Theta_X.$
Conversely, a mass function $m^X$ on $\Theta_X$ can be expressed in $\Theta_X\times \Theta_Y$  by moving each mass $m^X(B)$ to the cylindrical extension $B^{\uparrow XY} \doteq B \times \Omega_Y$ of $B$.
The \emph{vacuous extension} of $m^X$ onto $\Theta_X\times \Theta_Y$ is then: $m_X^{\uparrow XY}(A) \doteq m^X(B)$ if $A=B \times \Omega_Y$, 0 else.
The related BF is denoted by $Bel_X^{\uparrow XY}$.

\section{Likelihood-based inference} \label{sec:soa-inference-statistical-likelihood}
\label{sec:inference-likelihood}

Given a parametric model (\ref{eq:parametric-model}), we want to identify (or compute the support for) the parameter values which better describe the available data. 
An initial proposal for a likelihood-based belief function inference was made by Shafer \cite{Shafer76}, and immediately supported by Seidenfeld \cite{Seidenfeld78} and Wasserman \cite{Wasserman90}. Most recently, this approach has been endorsed by Denoeux \cite{Denoeux2010}, whereas its axiomatic foundations have been criticised by Moral \cite{Moral20141591}. Indeed, a major objection to likelihood-based inference is associated with the lack of commutativity with combination and conditioning: these aspects are not covered here for lack of space.

\subsection{Inference from traditional likelihood}

Consider the following requirements \cite{DENOEUX20141535}. Given the parametric model (\ref{eq:parametric-model}):
(i) the desired belief function $Bel_\Theta(\cdot | x)$ on the space $\Theta$ of parameter values should be based on the likelihood function $L(\theta;x) \doteq f(x | \theta)$ only (\emph{likelihood principle});
(ii) when a Bayesian prior $P_0$ on $\Theta$ is available, combining it with $Bel_\Theta(\cdot | x)$ using Dempster's rule should yield the Bayesian posterior: $Bel_\Theta(\cdot | x)\oplus P_0 = P(\cdot|x)$ (\emph{compatibility with Bayesian inference});
(iii) among all the belief functions which meet the previous two requirements, $Bel_\Theta(\cdot | x)$ should be the least committed.

These constraints lead to uniquely identify $Bel_\Theta(\cdot | x)$ as the consonant belief function
whose contour function $pl(\theta|x)$ is equal to the {normalised likelihood}:
\begin{equation} \label{eq:contour-function-likelihood-inference}
 pl(\theta | x)=\frac{L(\theta;x)}{\sup_{\theta'\in\Theta} L(\theta';x)}.
\end{equation}
The associated plausibility function is:
$
Pl_\Theta(A|x) = \sup_{\theta\in A} pl(\theta|x)
$ $\forall A \subseteq \Theta,
$ 
with as multivalued mapping $\Gamma_x : \Omega \rightarrow 2^\Theta$:
$
\Gamma_x(\omega) = \big \{ \theta\in \Theta \big | pl(\theta|x) \ge \omega \big \},
$
where $\Omega=[0,1]$ and the source probability there is the uniform one.

\begin{exam}\hspace{-1.5mm}\emph{\textbf{: Bernoulli sample}.} \label{exa:bernoulli-likelihood}
Let $x = (x_1,\ldots, x_n)$ consist of independent Bernoulli observations and $\theta\in\Theta=[0,1]$ be the probability of success. A \emph{Bernoulli trial} (or binomial trial) is a random experiment with exactly two possible outcomes, `success' and `failure', in which the probability of success is the same every time the experiment is conducted (repeated trials are independent). The contour function (\ref{eq:contour-function-likelihood-inference}) of the belief function obtained by likelihood-based inference is:
$
pl(\theta|x)=\frac{\theta^y(1-\theta)^{n-y}}{\hat{\theta}^y(1-\hat{\theta})^{n-y}},
$
with $y=\sum_{i=1}^n x_i$ and $\hat{\theta}$ the MLE.
\end{exam}

Denoeux \cite{DENOEUX20141535} argued that the method can be extended to handle low-quality data (i.e., observations that are only partially relevant to the population of interest). 
An extension of the Expectation-Maximization (EM) algorithm 
can also be proposed which maximises a generalised likelihood function able to handle uncertain data \cite{Denoeux2010}. \cite{1257308} proposed an iterative algorithm to estimate belief functions 
in this generalised EM framework.

\subsection{Belief likelihood function}

The previous framework simply takes the notion of likelihood as a given, and constructs belief functions from an input likelihood function. However, as recently shown in \cite{cuzzolin2020geometry}, there is no reason why we should not properly define a \emph{belief likelihood function}, mapping a sample observation $x \in \X$ to a real number,
rather than use the conventional likelihood to construct belief measures.
It is natural to define it as a family of BFs on $\X$, $Bel_\X(.|\theta)$, parameterised by $\theta \in \Theta$. 

In particular,
the belief likelihood function $Bel_{\X_1 \times \cdots \times \X_n} : 2^{\X_1 \times \cdots \times \X_n} \rightarrow [0,1]$ of a series of trials, $x_1,...,x_n$ is:
\[
Bel_{\X_1 \times \cdots \times \X_n}(A|\theta) \doteq Bel_{\X_1}^{\uparrow \times_i \X_i} \odot \cdots \odot Bel_{\X_n}^{\uparrow \times_i \X_i} (A|\theta),
\]
where $ Bel_{\X_j}^{\uparrow \times_i \X_i}$ is the vacuous extension of $Bel_{\X_j}$ to the Cartesian product $\X_1 \times \cdots \times \X_n$ where the observed tuples live, and $\odot$ is an arbitrary combination rule.
\\
Such a parameterised family is the input to the Generalised Bayesian Theorem \cite{smets93belief}. It takes values on \emph{sets} of outcomes, $A \subset \X$, of which singleton outcomes are mere special cases, providing a {natural setting for computing likelihoods of set-valued observations}, in accordance with the random set philosophy.
In particular, when $A = \{ \vec{x} \} = \{(x_1,...,x_n)\}$ reduces to a `sharp' sample, one can define its \emph{lower likelihood} $\underline{L}(\vec{x}) \doteq Bel_{\X_1 \times \cdots \times \X_n}(\{(x_1,...,x_n)\}|\theta)$ and \emph{upper likelihood} $\overline{L}(\vec{x}) \doteq Pl_{\X_1 \times \cdots \times \X_n}(\{(x_1,...,x_n)\}|\theta)$. 
When applied to samples generated by series of independent trials, under a generalisation of  stochastic independence, belief likelihood functions factorise into simple products. 
The resulting lower and upper likelihoods can be easily computed for series of Bernoulli trials.
This allows us to formulate a \emph{generalised logistic regression} framework \cite{cuzzolin2018belief,cuzzolin18belief-maxent}, in which the mass values of individual trials are constrained to follow a logistic dependence on scalar parameters.

\section{(Robust) Bayesian inference} \label{sec:robust-bayesian-inference}

Wasserman \cite{Wasserman90prior} noted that the mathematical structure of belief functions makes them suitable for generating classes of prior distributions, to be used in robust Bayesian inference. 
In particular, the upper and lower bounds of the posterior probability of a (measurable) subset of the parameter space $\Theta$ may be calculated directly in terms of upper and lower expectations.

Indeed, if a prior cannot be accurately specified, we might consider a credal set (an `envelope', in Wasserman's terminology) $\Pi$ and update each probability there using Bayes' rule, obtaining a new envelope $\Pi_x$ conditioned on $x \in X$.
\\
Denote $P_*(A) \doteq \inf_{P \in \Pi} P(A)$ and $P^*(A) \doteq \sup_{P \in \Pi} P(A)$ the lower and upper bounds on the prior induced by the envelope $\Pi$, and define
$L_A(\theta) = L(\theta) I_A(\theta)$, where $L(\theta) = f(x|\theta)$ is the likelihood and $I_A$ is the indicator function on $A$.
\begin{proposition} \label{pro:wasserman}
Let $\Pi$ be the credal set associated with a belief function on $\Theta$ induced by a source probability space $(\Omega,\mathcal{F}(\Omega),P_\omega)$
via a multivalued mapping $\Gamma$.
If $L(\theta)$ is bounded, then for any $A \in \mathcal{F}(\Theta)$: $\inf_{P_x \in \Pi_x} P_x (A) = \frac{E_*(L_A)}{E_*(L_A) + E^*(L_{A^c})}$, $\sup_{P_x \in \Pi_x} P_x (A) = \frac{E^*(L_A)}{E^*(L_A) + E_*(L_{A^c})}$,
where $E^*(f) = \int f^*(\omega) P_\omega(d\omega)$, $E_*(f) = \int f_*(\omega) P_\omega(d\omega)$ are the \emph{upper and lower integrals} of any real function $f$ on $\Theta$, where $f^*(\omega) \doteq \sup_{\theta \in \Gamma(\omega)} f(\theta)$ and $f_*(\omega) \doteq \inf_{\theta \in \Gamma(\omega)} f(\theta)$.
\end{proposition}
The integral representation for probability measures consistent with a belief function given by Dempster \cite{dempster2008upper} 
can also be extended to infinite sets (\cite{Wasserman90prior}, Theorem 2.1).

\section{Fiducial inference} \label{sec:inference-fiducial}

\subsection{Dempster's auxiliary variable} \label{sec:soa-inference-statistical-dempster}
\label{sec:inference-dempster}

In Dempster's approach to inference, the parametric (sampling) model (\ref{eq:parametric-model}) is supplemented by an a-equation (\ref{eq:a-equation}). The latter defines a multi-valued mapping $\Gamma : U \rightarrow 2^{\X \times \Theta}$ as:
\begin{equation} \label{eq:fiducial-random-set}
\Gamma(u) = \{ (x,\theta) \in \X \times\Theta | x = a(\theta,u) \} \subset \X \times \Theta. 
\end{equation}
Under standard measurability conditions \cite{Nguyen78,Molchanov05}, the probability space $(U, \mathcal{F}(U),\mu)$ and the mapping $\Gamma$ induce a belief function \textbf{$Bel_{\Theta\times\X}$} on $\X \times \Theta$. 

Conditioning (by Dempster's rule) $Bel_{\Theta\times\X}$ on $\theta$ yields the desired (belief) sample distribution on $\X$, namely:
$Bel_\X (A|\theta) = \mu(\{u: a(\theta,u) \in A\})$, $A \subseteq \X$.
Conditioning it on $X=x$ yields instead a belief measure $Bel_\Theta(\cdot | x)$ on $\Theta$, where $M_x(u) \doteq \{\theta : x = a(\theta,u)\}$:
\begin{equation} \label{eq:dempster-posterior}
Bel_\Theta(B| x) = \frac{ \mu (\{ u: M_x(u)\subseteq B \}) }{ \mu (\{ u: M_x(u) \neq \emptyset \}) }, \quad B \subseteq \Theta.
\end{equation}

\begin{exam}\hspace{-1.5mm}\emph{\textbf{: Bernoulli sample}.} \label{exa:bernoulli-dempster}
In the same situation as Example \ref{exa:bernoulli-likelihood}, consider the sampling model: $X_i =1$ if $u_i \le \theta$, 0 otherwise, where $u = (u_1,\ldots,u_n)$ has pivotal measure the uniform one on $[0,1]^n$: $\mu=\mathcal{U}([0,1]^n)$. Having observed the number of successes $y=\sum_{i=1}^n x_i$, the belief function $Bel_\Theta(\cdot|x)$ is induced by a random closed interval: 
$[u_{(y)},u_{(y+1)}]$,  
where $u_{(i)}$ denotes the i-th order statistics from $u_1,\ldots,u_n$. Quantities such as $Bel_\Theta([a,b]|x)$ or $Pl_\Theta([a,b]|x)$ can then be readily calculated.
\end{exam}

The issue of constructing belief functions via Dempster's model was later explored in \cite{almond92fiducial}, with a focus on belief function models for Bernoulli and Poisson processes.

\subsection{Inferential models, weak and elastic belief} \label{sec:inference-im} \label{sec:inference-weak-beliefs}

The \emph{inferential model} (IM) inference approach \cite{martin2010,Zhang11weak}
builds on the fiducial idea of trying to accurately predict the value $u^*$ of the auxiliary variable before conditioning on $X = x$, as more information is available about the former thanks to $\mu$.
For better predicting $u^*$,  the authors adopt a so-called `{predictive random set}' which  `smears'  the pivot distribution. 

Namely, given $Bel_\Theta(.|x)$, Dempster's posterior on $\Theta$ given the observable $x$ (\ref{eq:dempster-posterior}), we say that another belief function $Bel^*$ on $\Theta$ specifies an `inferential model' there if: $Bel^*(A) \leq Bel_\Theta(A)$ $\forall A \subseteq \Theta$.
\emph{Weak belief} (WB) is a method for specifying a suitable belief function within an inferential model. 
Weak belief weakens (hence the name) the `continue to believe' assumption that $u^*$ can be predicted by taking draws $u \sim \mu$ from the pivotal measure $\mu$, by choosing a set-valued mapping $\mathcal{S} : U \rightarrow 2^{U}$ that satisfies $u \in \mathcal{S}(u)$. As this is a multivalued mapping from the domain of the auxiliary variable onto itself, the pivot distribution $\mu$ induces a belief function on (the power set of) its own domain.
\\
The quantity $\mathcal{S}$ is called a \emph{predictive random set}, allowing uncertainty on the pivot variable itself.
At this point we have three domains $\X, U$ and $\Theta$ and two multi-valued mappings, $\Gamma : U \rightarrow 2^{\X \times \Theta}$ and $\mathcal{S} : U \rightarrow 2^{U}$.
The two belief functions on $U$ and $\X \times \Theta$ can be extended to $U \times \X \times \Theta$. By combining them there, and marginalising over $U$, we obtain a belief function whose random set:
$\Gamma_\mathcal{S}(u) = \bigcup_{u' \in \mathcal{S}(u)} \Gamma(u')$
is, by construction, dominated by the original Dempster's posterior. I.e., it is a valid inferential model.

In \cite{martin2010} various methods for building the mapping $\mathcal{S}$ are discussed.
In \cite{Zhang11weak} the authors illustrate their weak belief approach under: (i) inference about a binomial proportion, and (ii) inference about the number of outliers $(\mu_i \neq 0)$ based on the observed data $x_1, \cdots ,x_n$ under the model $X_i \sim \mathcal{N}(\mu_i,1)$.
Further extensions were proposed by Martin \cite{martin2010} and Ermini \cite{leaf12-inference}, under the name of \emph{elastic belief}. 

\subsection{Statistical inference with hints}

A variation of the fiducial argument was proposed by 
\cite{km95book}, who considered \emph{functional models}:
\begin{equation} \label{eq:functional-model}
f: \Theta \times \Omega \rightarrow \X, \quad x = f(\theta,\omega),
\end{equation}
by which observations $x \in \X$ are generated from a parameter $\theta \in \Theta$ and a random element $\omega \in \Omega$ with probability measure $P:\Omega \rightarrow [0,1]$, which plays the role of the pivot variable.

The observation $x$ induces an event in $\Omega$, namely:
$v_x \doteq \{ \omega \in \Omega | \exists \theta \in \Theta : x = f(\theta,\omega) \}$, which
in a Bayesian setting leads to condition the prior probabilities $P(\omega)$ with respect to $v_x$, obtaining $P'(\omega) = P(\omega) / P(v_x)$.
Assuming $\omega \in v_x$ had generated the observation, the possible values for the parameter $\theta$ are:
$T_x(\omega) = \{ \theta \in \Theta | x = f(\theta,\omega) \}$.
Summarising, an observation $x$ in a functional model (\ref{eq:functional-model}) generates a structure $\mathcal{H}_x = (v_x, P',T_x,\Theta)$ which Kohlas and Monney call a \emph{hint}.
We can then assess any hypothesis $H \subseteq \Theta$ on the correct value of the parameter with respect to it.
The arguments \emph{for} the validity of $H$ are the elements of the set $u_x(H) = \{\omega \in v_x : T_x(\omega) \subseteq H\}$, with degree of belief $P'(u_x(H))$, those merely \emph{compatible} with $H$ are $v_x(H) = \{\omega \in v_x : T_x(\omega) \cap H \neq \emptyset\}$, with plausibility $P'(v_x(H))$.
As in Equation (\ref{eq:fiducial-random-set}), the model
(\ref{eq:functional-model}) can itself be represented by a hint:
$\mathcal{H}_f = (\Omega,P,\Gamma_f, \X \times \Theta)$,
where
$\Gamma_f(\omega) = \{ (x,\theta) \in \X \times \Theta | x = f(\theta,\omega) \}$,
whereas an observation $x$ can be represented by the hint:
$\mathcal{O}_x = (\{v_x\},P,\Gamma, \X \times \Theta)$,
where $P(\{v_x\}) =1$, and $\Gamma(v_x) = \{x\} \times \Theta$.
The two pieces of information can then be combined, yielding: $\mathcal{H}_f \oplus \mathcal{O}_x$. By marginalising on $\Theta$ we obtain the desired information on the value of $\theta$: it is easy to show that the result is $\mathcal{H}_x$.

\section{Frequentist inference} \label{sec:inference-frequentist}


\subsection{A frequentist theory of lower probability} \label{sec:walley}

In a work predating Walley's theory of imprecise probability \cite{walley91book,Walley82frequentist} attempted to formulate a frequentist theory for upper and lower probability. 

Formally, the problem is to estimate a lower probability
$\underline{P}$ from a series of observations $x_1,...,x_n$. $\underline{P}$ is called a \emph{lower probability} whenever $\underline{P}(A \cup B) \geq \underline{P}(A) + \underline{P}(B)$ for all $A \cap B = \emptyset$ (super-additivity). $\overline{P}$ is an \emph{upper probability} whenever $\overline{P}(A \cup B) \leq \overline{P}(A) + \overline{P}(B)$, $A \cap B = \emptyset$ holds instead. Clearly, belief functions are a special case of lower probabilities.
\\
Walley and Fine considered the following estimator for $\underline{P}$:
\[
\underline{r}_n(A) \doteq \min \{ r_j(A) : k(n) \leq j \leq n \}, \quad k(n)
\rightarrow \infty
\]
where $r_j(A)$ is the relative frequency of event $A$ after observing $x_j$, and $k: \mathbb{N} \rightarrow \mathbb{N}$ goes to $\infty$ whenever $n \rightarrow \infty$.
Let $\underline{P}^\infty$ be the infinite IID product of the lower probability $\underline{P}$, describing the repetitions $x_1,x_2,...$ (\cite{Walley82frequentist}, page 746). This estimation process succeeds, in the sense that:
\[
\lim_{n \rightarrow \infty} \underline{P}^\infty(G_{n,\delta}^c) / \underline{P}^\infty(G_{n,\delta}) = 0 \quad \forall \delta>0,
\]
where $G_{n,\delta}$ is the event $|\underline{r}_n(A) - \underline{P}(A)|<\delta$.
This result parallels Bernoulli's law of large numbers: the confidence that $\underline{r}_n(A) $ is close to $\underline{P}(A)$ grows with the sample's size\footnote{Note that the opposite view was supported in \cite{lemmers86confidence}.}.

\subsection{Dempster's (p,q,r) interpretation}

In more recent times, \cite{
DEMPSTER2008365} proposed a semantics for belief functions whereby every assertion $A$ is associated with a triple $(p,q,r)$ where $p = Bel(A)$ is the probability `for' the assertion, $q = Bel(A^c)$ is the probability `against' the assertion, and $r = 1 - p - q$ is the probability of `don’t know'. 
The methodology was applied to inference and prediction from Poisson counts, and the relation of DS theory to statistical significance testing elaborated, by introducing the concept of \emph{dull} null hypothesis (i.e., a hypothesis which assigns mass 1 to an interval of values). 
Poisson's a-probabilities and values for $p,q$ and $r$ were then derived.

\subsection{From confidence intervals}

Confidence intervals (\ref{eq:confidence-intervals}) can also be exploited to quantify beliefs about the realisation of a discrete random variable $X$ with unknown probability distribution $P_X$. In \cite{DENOEUX2006228} a solution which is less committed than $P_X$, and converges towards the latter in probability as the size of the sample tends to infinity, was proposed. Namely, each confidence interval can be thought of as a credal set, specifically a set of (feasible) probability intervals:
$\mathcal{P}(l,u) \doteq \{p: l(x) \leq p(x) \leq u(x), \forall x \in \Theta \}$.
A set of probability intervals is called \emph{feasible} \cite{decampos94} if and only if for each
$x \in \Theta$ and every value $v(x) \in [l(u), u(x)]$ there exists a probability distribution function $p:\Theta \rightarrow [0,1]$ for which $p(x) = v(x)$.
One can then obtain the lower and upper probabilities (see Section \ref{sec:walley}) associated with a confidence interval using the following simple formulas:
$\lpr(A)=\max \{ \sum_{x \in A} l(x), 1- \sum_{x \not\in A} u(x) \}$, 
$\lpr(A)=\min \{ \sum_{x \in A} u(x), 1- \sum_{x \not\in A} l(x) \}$.

\subsection{A theory of confidence structures} 

A recent work \cite{BALCH20121003} attempted to reconcile belief function theory with frequentist statistics at a more fundamental level, using conditional random sets. Once again, given a statistical model 
(\ref{eq:parametric-model}) 
we wish to construct a BF $Bel_\Theta(.|x) : 2^\Theta \rightarrow [0,1]$.
An \emph{observation-conditional random set} on the parameter $\theta$ is a mapping $\Gamma: \Omega \times \X \rightarrow \Theta, \Gamma(\omega,x) \subset \Theta$ which depends also on the observable $x$. The resulting belief value is:
$Bel_\Theta(A|x) = P_\omega(\{\omega \in \Omega : \Gamma(\omega,x) \subseteq A\})$
for all $x \in \X$, $A \subset \Theta$. The problem then reduces to that of constructing an appropriate observation-conditional random set from the given statistical model.

A \emph{confidence set} is a set-estimator for the parameter, a function of the observable
$C :  \X \rightarrow \mathcal{P}(\Theta)$, $C(x) \subset \Theta$
designed to cover the true parameter value with a specified regularity. Confidence intervals are the most well-known type of confidence set. The \emph{Neyman-Pearson confidence} associated with a confidence set is then the frequentist probability of drawing $x$ such that $C(x)$ covers the true parameter value.
A \emph{confidence structure} is an observation-conditional random set whose source probabilities $P_\omega$ are commensurate with Neyman-Pearson's confidence:
$P_\X ( \{x \in \X : \theta \in \bigcup_{\omega \in A} \Gamma(\omega,x) \} | \theta ) \geq P_\omega(A)$ for all $A \subset \Omega, \theta \in \Theta$.
In other words, a confidence structure is a means for constructing confidence sets of the form:
$C(A;x) =  \bigcup_{\omega \in A} \Gamma(\omega,x)$, $A \subset \Omega,$
conditioned on the coverage probability of $C(A;x)$ being greater than or equal to the source probability, $P_\omega(A)$.
Since every $B \subset \Theta$ can be put in the form $B=C(A;x)$ for some $A \subset \Omega$, $x \in \mathcal{X}$, all sets of parameter values are confidence intervals, with their belief values  measuring how much confidence is associated with each such set. 
In \cite{BALCH20121003} various (inference) methods for constructing confidence structures from confidence distributions (Section 3.1), pivots (3.2) and p-values (Section 3.3) were illustrated,
and contrasted with other statistical inference methods.


\section{Other statistical approaches} \label{sec:soa-inference-statistical-others}

Various other scholars contributed to the study of statistical inference with belief functions \cite{acker00belief}.

A belief function generalisation of Gibbs ensembles was proposed by \cite{Kong1988Gibbs}.
\cite{gillett2000attribute} discussed the integration of statistical evidence from attribute sampling with non-statistical evidence within the belief function framework. They also showed how to determine the sample size in attribute sampling to obtain a desired level of belief that the true attribute occurrence rate of the population lies in a given interval, and what level of belief is obtained for a specified interval given the sample result.
\cite{Edlefsen2009} presented an approach to estimating limits from Poisson counting data with nuisance parameters, by deriving a posterior belief function for the `Banff upper limits challenge' three-Poisson model.
Novelty detection, i.e., the problem of testing whether an observation may be deemed to correspond to a given model, was also studied in the belief framework \cite{Aregui2006NoveltyDI}.


\section{Shafer's reflections}

Shafer himself had analysed back in 1982 three ways of doing statistical inference with belief functions \cite{Shafer82},
according to the nature of the evidence inducing the parametric model. When
distinct, independent observations for each $f(.|\theta)$ (e.g. symptoms $x$ arising from a medical condition $\theta$) are available, he argued, Smets' \emph{conditional embedding} \cite{smets93belief} should be adopted, which consists of extending each $P_\theta = f(.|\theta)$ on $\X \times \Theta$ and combining them there by Dempster's rule. 
The case when the parametric model is induced by error distributions should be treated using fiducial inference (Section \ref{sec:inference-fiducial}).

More interesting is the case in which  no evidence is there, except the convinction that the phenomenon is random.
This case can draw inspiration from the Bayesian treatment there, in which the random variable $X$ is thought of as one of a sequence $\mathbf{X} = (X_1,X_2,...)$ of unknown quantities taking values in $\X = \{x_1,...,x_k\}$, and the Bayesian beliefs on $\mathbf{X}$ are expressed by a countably additive and symmetric probability distribution $P$. Additivity and symmetry imply that limiting frequencies exist $P(\lim_{n \rightarrow \infty} f(x,n) \;\text{exists}) = 1$, where $f(x,n)$ is the proportion of the quantities $X_1,...,X_n$ that equal $x \in \X$, and that the probability distribution $P_\theta, \theta \in \Theta$ on the parameter can be recovered by conditioning $P$ on these limiting frequencies.
Shafer thus argued in favour of constructing a belief function for $\mathbf{X} = (X_1,X_2,...)$ such that: (1) $Bel(\lim_{n \rightarrow \infty} f(x,n) \; \text{exists} ) =1$ for all $x \in \X$; (2)
$Bel (X_1=x_1, ..., X_n = x_n | \lim_{n \rightarrow \infty} [f(1,n),...,f(k,n)] = \theta )$ 
is equal to $P_\theta(x_1), ..., P_\theta(x_n)$ for all $x_1,...,x_n \in \X$;
(3) $Bel(\lim_{n \rightarrow \infty} [f(1,n),...,f(k,n)] \in A) = 0$ $\forall A \subsetneq \Theta$.
Such belief functions exist, and can be constructed. Interestingly, in the case $\X = \{0,1\}$ the solution is unique and its marginal for $(\theta, X_1)$ is a BF on $\X \times \Theta$ which corresponds to Dempster's generalised Bayesian solution. In the general case $k>2$, unfortunately, the construction is less direct. 
Additionally, no proofs of the above arguments were given in \cite{Shafer82}.

\section{Discussion}
More recent developments allow us to have a clearer picture of the topic. 
Wasserman's credal approach, whereas sound and rooted on the standard Bayesian formalism, has not gained much ground. The main reason is that, being designed to generic credal sets, it is simply not specifically tailored for envelopes associated with belief functions. Moreover, the credal interpretation is not compatible with Dempster-Shafer theory as a whole, with its specific combination rule.
A lot of effort, led by Dempster, has been directed in the fiducial sense.
Dempster's model indeed allows us to quantify the uncertainty on $\Theta$ without having to specify a prior distribution on $\Theta$, and is compatible with Bayesian inference in the sense defined in Section \ref{sec:inference-likelihood}.
However, it often leads to cumbersome calculations which require the use of Monte-Carlo simulations. More fundamentally, the fiducial argument requires an a-equation 
and a pivot variable 
which are not observable, and not uniquely determined for a given statistical model. 
Statistical inference with hints, being roughly equivalent to Dempster's approach, is open to the same criticism.
The issue is even more acute in the weak belief setting, which requires an additional structure in the form of a predictive random set, itself a rather complex construct with no link to observables.


Most frequentist efforts produce lower and upper probabilities, rather than random sets \cite{cuzzolin2023reasoning} (belief functions). This Author believes the notion of extending statistical hypothesis testing to the belief framework (as in Dempster's dull hypothesis), on the other hand, quite interesting. Confidence structure theory also seems promising, but suffers from the same model overfit issue as the fiducial approach, for observation-conditional random sets are a design choice. It would interesting to investigate its extension to Bayesian \emph{credible intervals}\footnote{\scriptsize \url{https://www.coursera.org/learn/bayesian/lecture/hWn0t/credible-intervals}}, as well. 

The likelihood-based approach, instead, despite having been supported by many, has a clear limitation in the sense that it produces only consonant belief functions (which are in 1-1 correspondence with possibility distributions). Thus, it does not exploit the full expressive power of belief theory, but should probably be considered an approach to possibilistic, rather than belief function, inference.
An interesting discussion of the connection between belief functions and the likelihood principle can be found in \cite{Aickin2000}. There it is argued that expressing a belief value about a parameter in a statistical model is not consistent with the likelihood principle, as the ordering of the operations (combining belief, combining likelihood) makes a difference.
The belief likelihood idea, in opposition, does not require extra hidden variables and structures, and delivers fully fledged belief functions.

\section{Conclusions}

In this survey we appreciated how the vast majority of approaches to belief function inference from statistical data can be grouped according to what mainstream probability interpretation they generalise: likelihood-based, Bayesian, fiducial, or frequentist. We recalled all related basic notions, provided some simple examples, and discussed their relative strengths and limitation, which led us to identify some methodologies as particularly promising. 

Due to lack of space, related issues such as belief function inference from qualitative or partial data, or in the machine learning \cite{cuzzolin13fusion,manchingal2025epistemic-b,cuzzolin2024epistemic,Manchingal2025AUE,wang2025review,caprio2025credal} or statistical learning theory \cite{cuzzolin2024generalising} context are not covered here. Indeed, 
as argued in \cite{Shafer82}, the strength of belief calculus is really about allowing inference under partial knowledge or ignorance, when simple parametric models are not available.

\bibliographystyle{named}
\bibliography{bibliography,new_citations}

\end{document}